\documentclass[11pt,reqno]{amsart}

\usepackage[centertags]{amsmath}
\usepackage{amsfonts}
\usepackage{amssymb}
\usepackage{amsthm}
\usepackage{newlfont}
\usepackage{url}
\usepackage{setspace}

\usepackage{graphicx}
\usepackage[usenames]{color}
\usepackage{caption}
\usepackage{subcaption}
\usepackage{enumitem}
\usepackage{wrapfig}
\usepackage{verbatim}

\theoremstyle{plain}
\newtheorem{thm}{Theorem}

\theoremstyle{definition}
\newtheorem{defn}{Definition}
\theoremstyle{remark}

\numberwithin{thm}{section}
\numberwithin{defn}{section}

\newcommand{\defeq}{\mathrel{\mathop:}=}

\def\D{\mathbb D}

\def\R{\mathbb R}
\def\Z{\mathbb Z}

\begin{document}
\title[]{Link homotopic but not isotopic}
{\hfil }
\author[]{Bakul Sathaye}

\pagestyle{headings}
\begin{abstract} 
Given an $m$-component link $L$ in $S^3$ ($m \ge 2$), we construct a family of links which are link homotopic, but not link isotopic, to $L$. Every proper sublink of such a link is link isotopic to the corresponding sublink of $L$. Moreover, if $L$ is an unlink then there exist links that in addition to the above properties have all Milnor invariants zero. 
\end{abstract}


\maketitle

\section{Introduction}

An $n$-component link $L$ in $S^3$ is a collection of piecewise linear maps $(l_1, \ldots, l_n) : S^1 \rightarrow S^3$, where the images $l_1 (S^1), \ldots , l_n(S^1)$ are pairwise disjoint. A link with one component is a knot.
 
Two links in $S^3$, $L_1$ and $L_2$,  are said to be \textit{isotopic} if there is an orientation preserving homeomorphism $h: S^3 \rightarrow S^3$ such that $h(L_1) = L_2$ and $h$ is isotopic to the identity map. 
 
The notion of link homotopy was introduced by Milnor in \cite{Mil}. Two links $L$ and $L'$ are said to be \textit{link homotopic} if there exist homotopies $h_{i,t}$, between the maps $l_i$ and the maps $l_i'$ so that the sets $h_{1,t}(S^1), \ldots ,h_{n,t}(S^1)$ are disjoint for each value of $t$. In particular, a link is said to be link homotopically trivial if it is link homotopic to the unlink. Notice that this equivalence allows self-crossings, that is, crossing changes involving two strands of the same component.

The question that arises now is: how different are these two notions of link equivalence? 
From the definition it is clear that a link isotopy is a link homotopy as well. Link homotopy allows self-crossings, which dramatically simplifies the equivalence between links. Under link homotopy, all knots become trivial knots, despite the fact that they are non-trivial under isotopy. Link homotopy in some sense measures the linking between different components of a link. So the question reduces to whether two link homotopic links are also isotopic. We have the example of the Whitehead link, which is a non-trivial link that becomes trivial under link homotopy. This gives us an example of a 2-component link that is link homotopic, but not isotopic, to the unlink. But there are no known examples of such links with more than two components. 

The main theorem in this article provides a method to generate examples of links with more than two components that are link homotopic but not isotopic to any given link. In particular, we prove the following.

\begin{thm}
\label{mainthm}
\textit{Let $L$ be any $m$-component link in $S^3$ for $m\ge 2$. There exists a family of $m$-component links, $\mathcal L_m$, with the following properties: }
\vskip -2pt
\begin{enumerate}
\itemsep 2pt
\item \textit{Each $\mathcal L_m$ is link homotopic to $L$ }
\item \textit{Each $\mathcal L_m$ is not link isotopic to $L$}
\item \textit{Every proper sublink of $\mathcal L_m$ is link isotopic to the corresponding sublink of $L$. In particular, if $L$ is an unlink then every proper sublink of $\mathcal L_m$ is isotopic to the unlink.} 
\end{enumerate} 
\end{thm}

The proof of this theorem constitutes Sections 4 and 5. In the last section we show how to generate examples of links with an additional property about Milnor invariants. 
We prove the following.

\newpage
\begin{thm}
\label{milnorinv}
There exist links in $S^3$ with the following properties. 
\begin{enumerate}
\item \textit{It is link homotopic to the unlink.}
\item \textit{It is not link isotopic to the unlink.}
\item Every proper sublink is link isotopic to the unlink. 
\item \textit{All Milnor's $\mu$-invariants, $\overline{\mu}(i_1 \ldots i_s)$, are $0$. }
\end{enumerate}
\end{thm}


\vspace{4mm}

\noindent \textbf{Acknowledgements.}
I am grateful to my advisor, Jean-Fran\c{c}ois Lafont, for posing the question addressed in this paper and for his support and guidance through the process. I would like to thank J.-B. Meilhan for his inputs. 

\section{Jones polynomial} 

The Jones polynomial gives an invariant of an oriented link, that assigns to each link a Laurent polynomial in the variable $t^{1/2}$ with integer coefficients. If two links are link isotopic then they have the same Jones polynomial.  
 
The Jones polynomial is a function 
$$V : \{ \mathrm{Oriented}\ \mathrm{links}\ \mathrm{in}\ S^3 \} \rightarrow \mathbb Z [t^{-1/2}, t^{1/2} ]$$ 
 with the following properties. 
\begin{enumerate} 
\item $V(\mathrm{unknot}) = 1$.
\item $V( L \cup \bigcirc) = (-t^{-1/2} - t^{1/2}) V(L)$.
\item $V$ satisfies the Skein relation: whenever three oriented links $L_+, L_-$ and $L_0$ are the same, except in the neighbourhood of a point where they are as shown in Figure \ref{Skein relation}, then 
$$t^{-1}V(L_+) - tV(L_-) + (t^{-1/2} - t^{1/2}) V(L_0) = 0.$$ 
 
\begin{figure}[h!]
\includegraphics[scale=0.8]{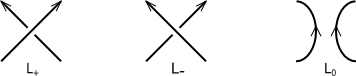}
\caption{Skein relation }
\label{Skein relation}
\end{figure}
 \end{enumerate} 
 
The Jones polynomial is discussed in detail in the book \cite{Lic}.

\vspace{2mm}

\section{Whitehead link} 

\begin{wrapfigure}{l}{0.39\textwidth}
\vspace{-2mm}
\centering
\includegraphics[scale=0.7]{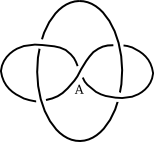}
\caption{Whitehead link}
\label{Whitehead link}
\vspace{-5mm}
\end{wrapfigure}

The Whitehead link is given in Figure \ref{Whitehead link}. It can be seen to be link homotopic to the unlink by making a crossing change at the central crossing, $A$. However, its Jones polynomial is computed to be $t^{-3/2}(-1 + t -2t^2 + t^3 - 2t^4 + t^5)$, while the Jones polynomial of the unlink is $(-t^{-1/2}-t^{1/2})$.

 
This shows that the Whitehead link is not link isotopic to the unlink. Moreover, removing any component from the Whitehead link leaves us with the unknot. Thus, it satisfies the conclusion of Theorem \ref{mainthm} for $m=2$. 
 
\newpage
\section{Unlink}
 
In this section, we construct examples of links as in Theorem \ref{mainthm} where $L$ is an unlink. 

\subsection{Construction} 
\begin{figure}[h!]
\includegraphics[scale=0.8]{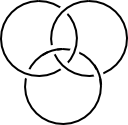}
\caption{Borromean rings}
\label{Borromean}
\end{figure}

\begin{figure}[h!]
\includegraphics[scale=0.7]{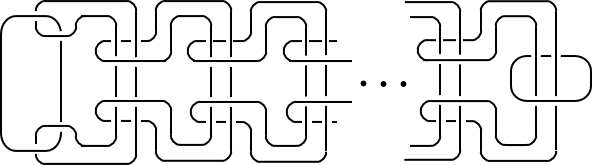}
 \caption{Brunnian link with $m$ components, $B_m$}
 \label{Brunnian-n}
\end{figure}

A \textit{Brunnian link} is a non-trivial link such that each sublink is link isotopic to an unlink. An example of a Brunnian link with 3 components are the Borromean rings (Figure \ref{Borromean}).  For $m\ge 3$, a family of Brunnian links, $B_m$, was given by Milnor \cite{Mil} as shown in Figure \ref{Brunnian-n}. Notice that these links satisfy properties (2) and (3) of our theorem, but not property (1).

\begin{figure}[h!]
 \begin{subfigure}[b]{0.3\textwidth}
        \includegraphics[width=\textwidth]{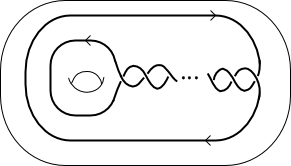}
        \caption{$n$ negative twists}
        \label{fig:gull}
    \end{subfigure}
    \qquad \qquad \qquad 
      \begin{subfigure}[b]{0.3\textwidth}
        \includegraphics[width=\textwidth]{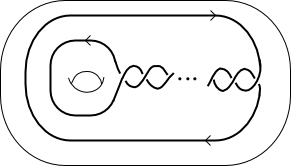}
        \caption{$n$ positive twists}
        \label{fig:tiger}
    \end{subfigure}
       \caption{Whitehead $n$ twists, $\mathcal W_n$}\label{fig:animals}
       \label{Whitehead double}
 \end{figure}

\begin{defn} [Whitehead double]
\label{Whiteheaddouble}
 Let $K$ be a component of a link $L$ in $S^3$, regarded as $h(\{ 0 \} \times S^1)$ for some embedding $h : D^2 \times S^1 \rightarrow S^3 \setminus (L\setminus K)$, such that $K$ and $h((0,1) \times S^1)$ have linking number zero. Let $n$ be a non-zero integer. Consider in the solid torus $T = D^2 \times S^1$ the knot $\mathcal W_n$ with $n$ crossings, as depicted in Figure \ref{Whitehead double}. The knot $h(\mathcal W_n)$ is called the Whitehead $n$-double of $K$, and is denoted by $W_n(K)$.  
\end{defn} 

If the components of the link $L$ are numbered, $L = K_1 \cup \ldots \cup K_m$, then the link obtained by considering the Whitehead $n$-double on the $i$-th component of $L$ will be denoted by $W_n^i(L)$. We are now ready to prove our theorem.


\subsection{Case $m = 3$}
\label{m=3}
\begin{figure}[h!]
\includegraphics[scale=0.7]{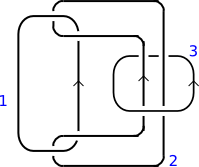}
\caption{Brunnian link with 3 components, $B_3$}
\label{B_3}
\end{figure}

\begin{figure}[h!]
 \includegraphics[scale=0.7]{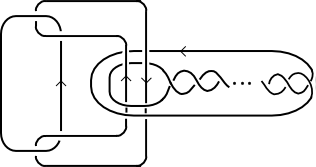}
 \caption{Whitehead double of $B_3$}
 \label{W_n(B_3)}
 \end{figure}

We first consider the case of 3-component links. Let $B_3$ be the Brunnian link with 3 components, with its components numbered as shown in Figure \ref{B_3}.
Consider the Whitehead $n$-double (with negative twists) of the third component of $B_3$ to obtain the new link, $W_n^3(B_3)$. To simplify notation, we will denote this by $W_n(B_3)$ hereafter (Figure \ref{W_n(B_3)}). 
 
Observe that when $n$ is even, this link is link homotopic to the unlink. Even a single crossing change in $W_n(B_3)$ unwinds the link to give an unlink. Note that $W_n(B_3)$ satisfies properties (1) and (3) of our theorem.
 
We claim that there is an integer $n>0$ for which $W_n(B_3)$ is not link isotopic to the unlink, that is, it satisfies property (2). In fact, such an $n$ occurs in every pair of consecutive even integers. We show this by showing that the Jones polynomial of $W_n(B_3)$ for some $n$ is not the same as that of the unlink. 

\begin{figure}[h!]
\includegraphics[width=\textwidth]{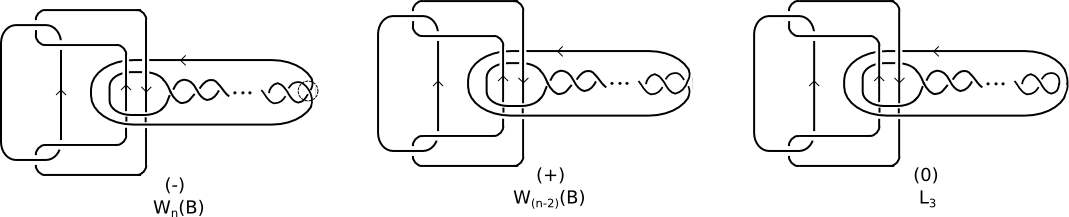}
\caption{}
\label{Skein for W_n}
 \end{figure}

Consider the skein relation at the crossing indicated in Figure \ref{Skein for W_n}.

\begin{equation}
t^{-1}V(W_{n-2}(B_3)) - tV(W_n(B_3)) + (t^{-1/2} - t^{1/2}) V(L_3) = 0.\label{skein}\end{equation}
 
If our claim is false, then both $W_{n-2}(B_3)$ and $W_n(B_3)$ must be link isotopic to the unlink and must have the same Jones polynomial as the unlink with 3 components, namely, $(-t^{-1/2} - t^{1/2})^2$.  

In that case, the above equation becomes 
\begin{equation}
t^{-1}(-t^{-1/2} - t^{1/2})^2- t(-t^{-1/2} - t^{1/2})^2 + (t^{-1/2} - t^{1/2}) V(L_3) = 0.
\end{equation}
Therefore, \begin{equation} (t^{-1/2} - t^{1/2}) V(L_3) = (t-t^{-1})(-t^{-1/2} - t^{1/2})^2. \end{equation}
That is, \begin{equation} V(L_3) = (-t^{-1/2} - t^{1/2})^3. \label{VL_3}\end{equation}
 
\begin{figure}[h!]
\includegraphics[scale=0.8]{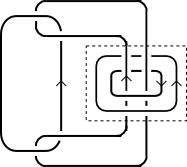}
\caption{}
\label{Boxed L_3}
 \end{figure}

\vskip 15pt 
Now we compute the Jones polynomial for $L_3$ through multiple applications of the skein relations.  
Consider the link $L_3$ and focus on the part of the link indicated in Figure \ref{Boxed L_3}.

We apply the skein relation at the crossings indicated to get the following equations.\\

\begin{center}
\includegraphics[scale=0.5]{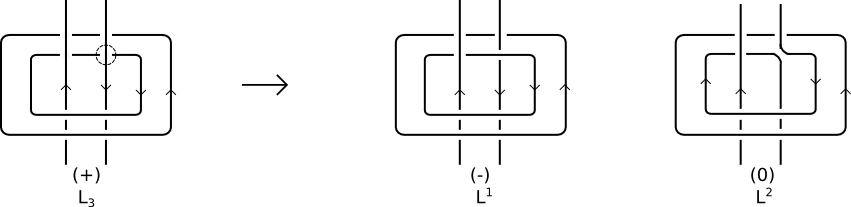}
\end{center}
\begin{equation} t^{-1} V(L_3) = t V(L^1) - (t^{-1/2} - t^{1/2}) V(L^2) \end{equation}
\hspace{1mm}
\begin{center}
\includegraphics[scale=0.5]{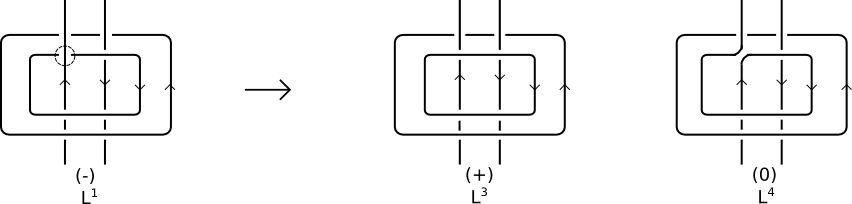}
\end{center}
\begin{equation} t V(L^1) = t^{-1} V(L^3) + (t^{-1/2} - t^{1/2}) V(L^4)\end{equation}
\hspace{1mm}
\begin{center}
\includegraphics[scale=0.5]{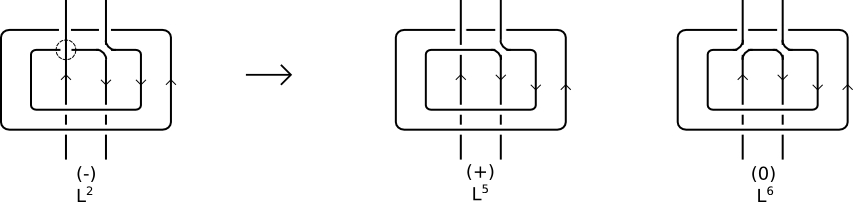}
\end{center}
\begin{equation} t V(L^2) = t^{-1} V(L^5) + (t^{-1/2} - t^{1/2}) V(L^6)\end{equation}

\begin{figure}[h!]
\begin{subfigure}[b]{0.3\textwidth}
        \includegraphics[scale=0.5]{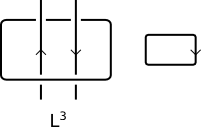}
         \end{subfigure}
    ~     \begin{subfigure}[b]{0.22\textwidth}
        \includegraphics[scale=0.5]{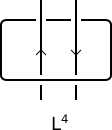}
    \end{subfigure}
     \begin{subfigure}[b]{0.2\textwidth}
        \includegraphics[scale=0.5]{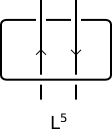}
         \end{subfigure}
     \begin{subfigure}[b]{0.25\textwidth}
        \includegraphics[scale=0.5]{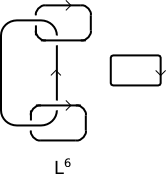}
     \end{subfigure}
      \caption{}
      \label{L^}
 \end{figure}

Under isotopy the links $L^3, L^4, L^5$ and $L^6$ can be drawn as shown in Figure \ref{L^}.

Now define the link $A$ as shown.

\begin{figure}[h!]
\includegraphics[scale=0.7]{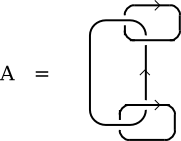}
\label{A}
\end{figure}

Notice that the links $L^3, L^4, L^5$ and $L^6$ can be written, up to isotopy, using the links $A$ and the Brunnian link, $B_3$. Indeed, we have
$L^3 = B_3 \cup \bigcirc $,
$ L^4 =  B_3$,
$L^5 =  B_3$,
$L^6  =  A \cup \bigcirc$.

So now our computation is reduced to finding the Jones polynomial for $B_3$ and $A$. After applying further Skein relations, we can calculate:

$$V(B_3) = - t^3 + 3t^2 - 2t + 4 - 2t^{-1} + 3t^{-2} - t^{-3} $$
$$V(A) = t^2 + 2 + t^{-2}.$$
Using the Skein relations above, this allows us the calculate:
\begin{equation}
V(L_3) = t^{9/2} - 2t^{7/2} + t^{5/2} - 2t^{3/2} - 2t^{1/2} - 2t^{-1/2} - 2t^{-3/2} + t^{-5/2} - 2t^{-7/2} + t^{-9/2}.
\label{L_3}
\end{equation}
However, this is different from the expression we obtained in (\ref{VL_3}) for $V(L_3)$, proving our assumption to be incorrect. Hence, both $W_n$ and $W_{n-2}$ cannot be unlinks. This shows that at least one of any pair of consecutive even numbers gives a link with properties as required in Theorem \ref{mainthm}. We have thus proved our theorem for links with $m=3$ components. 


\subsection{Case $m > 3$} 
\begin{figure}[h!]
\includegraphics[scale=0.7]{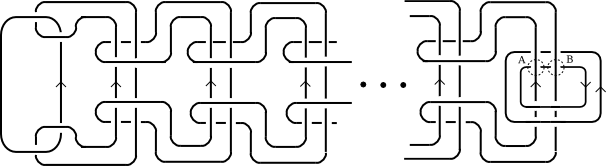}
\caption{}
\label{L_m}
\end{figure}

Now consider the case when the number of components is greater than 3. Consider the $m$-component Brunnian link, $B_m$ $( m >3)$, and the Whitehead $n$-double of the $m$-th component to obtain the link $W_n(B_m)$. Using the Skein relations, similar to the ones in Figure \ref{Skein for W_n}, we obtain the link $L_m$ as shown in Figure \ref{L_m}. 

If $W_n(B_m)$ and $W_{n-2}(B_m)$ were both isotopic to the unlink, then computations similar to equations $(\ref{skein}) - (\ref{VL_3})$ show that the Jones polynomial for $L_m$ will be 
\begin{equation} V(L_m) = (-t^{-1/2} - t^{1/2})^m. \label{VL_m}\end{equation}

On the other hand, we compute the Jones polynomial of $L_m$ by application of Skein relations at the crossing (A) and then at (B). This reduces the link $L_m$ to the links $L^1_m, \ldots ,L^6_m$, similar to the ones in equations $(5) - (7)$. Here
$L^3_m = B_m \cup \bigcirc,\ L^4_m = B_m,\ L^5_m = B_m$ and $L^6_m = L_{m-1} \cup \bigcirc.$

Computations similar to the ones in the case $m = 3$ give us the following recurrence relation for the Jones polynomial for $L_m$, for every $m>3$:
\begin{equation}
V(L_m) =   V(L_{m-1}) (t^{5/2} -t^{3/2} -t^{-3/2} +t^{-5/2}) + (-t^{1/2}-t^{-1/2})^m (t^2 - 2t+ 3 - 2 t^{-1} + t^{-2} )
    \label{recrelation}
\end{equation}

 

We claim that deg$(V(L_m)) >m/2$ for every $m\ge 3$. From equation (\ref{L_3}), we know that deg$\displaystyle (V(L_3)) = \frac{9}{2} > \frac{3}{2}$, hence the result is true for $m=3$. 

Assume that this is true for $m-1$. From the relation \eqref{recrelation} we see that the highest power of $t$ in $V(L_m)$ comes from one of the terms, $t^{5/2} V(L_{m-1})$ and $(-t^{1/2}-t^{-1/2})^mt^2$. By the induction hypothesis, we know that deg$\displaystyle (V(L_{m-1})) > \frac{m-1}{2}$. Hence, deg $\displaystyle (t^{5/2}V(L_{m-1})) > \frac{5}{2} + \frac{m-1}{2}= \frac{m+4}{2}$. This shows that degree of $V(L_m)$ must be the degree of the term $t^{5/2}V(L_{m-1})$, that is, deg$(V(L_{m}))  =$ deg$\displaystyle (V(L_{m-1})) + \frac{5}{2}$. Also, deg$\displaystyle (V(L_m)) > \frac{m+4}{2} > \frac{m}{2}$, which proves our claim.\\

It also gives us a recurrence relation for the degree of $V(L_m)$, namely, 
\begin{center}deg$(V(L_m)) =$ deg$\displaystyle (V(L_{m-1})) + \frac{5}{2}$.
\end{center} 
This in turn gives us that, for $m \ge 3$, deg$\displaystyle (V(L_m)) = \frac{5}{2}(m-3) + \frac{9}{2}$.  
 
This shows that \begin{center}degree of $V(L_m) >$ degree of $(-t^{1/2}- t^{-1/2})^m$, for $m\ge 3$,\end{center} showing that the expression obtained in \eqref{VL_m} must be different from the one in the recurrence relation \eqref{recrelation}. 

Thus one of the links, $W_n(B_m)$ or $W_{n-2}(B_m)$, must be non-isotopic to the unlink. This completes the proof of Theorem \ref{mainthm} for $m \ge 3$.  
 
\noindent \textbf{Remark.} 
J.-B. Meilhan pointed out that a result of Meilhan-Yasuhara (Theorem 1, \cite{MY}) can also be used to show that $W_n(B_m)$ are non-trivial links for every $n>0$. Their result shows that these links have non-zero Milnor isotopy invariants \cite{Mil2}, showing that they are link isotopically non-trivial. And hence they also fail to satisfy property (4) of Theorem \ref{milnorinv}. In Section \ref{milnor}, we give examples satisfying properties $(1) - (4).$

\vspace{3mm}
\subsection{Non-isotopic links}

Here we point out that the method described in this section is useful to further show that each of the links, $W_n(B_3)$, are in fact non-isotopic links. To see this, observe that one can obtain the Jones polynomial of $W_n(B_3)$  as a function of $W_{n-2}(B_3)$, using the Skein relation in equation (\ref{skein}). We will do that here.

Define the function $\Phi$ as follows.  $$\Phi: \Z [t^{-1/2}, t^{1/2}] \rightarrow \Z [t^{-1/2}, t^{1/2}] $$
$$\Phi(p) = t^{-2} p + t^{-1} (t^{-1/2} - t^{1/2}) V(L_3),$$
 where $V(L_3)$ is the Jones polynomial of the link $L_3$ as given in equation (\ref{L_3}). 
Then $V(W_n(B_3)) = \Phi(V(W_{n-2}(B_3)))$, for every $n \ge 2$. Hence we only need to analyse the function $\Phi$ in order to understand the Jones polynomials of the Whitehead doubles. 

One can see that the degree of $\Phi(p)$ depends on that of the polynomial $p$. Let us denote by deg$_L (p)$ the degree of the lowest degree term of $p$. Then, 
\begin{center} deg$_L(\Phi(p)) = $ min \{deg$_L (p) - 2, $ deg$_L (V(L_3)) - 3/2 $\}. \end{center} 
We know from our previous computations that deg$_L (V(L_3)) = -9/2$. Hence,
\begin{equation} \mathrm{deg}_L(\Phi(p)) = \mathrm{ min \{deg}_L (p) - 2, -6 \}.\label{reln}\end{equation} 

For $n=2$, we have that $V(W_2(B_3)) = \Phi(V(W_0(B_3)))$. Since $W_0(B_3)$ is nothing but the unlink with 3 components, $V(W_0(B_3)) = (-t^{-1/2} - t^{1/2})^2$, and deg$_L(V(W_0(B_3))) = -1$. Hence deg$_L(V(W_2(B_3))) =$ deg$_L(\Phi(V(W_0(B_3)))) = -6$.

Continuing in a similar fashion, we get that deg$_L(V(W_4(B_3))) = -8,$ deg$_L(V(W_6(B_3))) = -10$ and so on. 

We can do similar computations for odd $n$. Using skein relations, we get $V(W_1(B_3)) = - t^4 + 2t^3 - t^2 + 2t + 1 + t^{-2} - t^{-3} + 2t^{-4} - t^{-5}$. So, deg$_L(V(W_1(B_3))) = -5.$ Using the above relation \ref{reln} we get deg$_L(V(W_3(B_3))) =$ deg$_L(\Phi(V(W_1(B_3)))) = -7$, deg$_L(V(W_5(B_3))) = -9$ and so on. 

Together, we get that deg$_L(V(W_{n}(B_3))) = -n-4$ for $n\ge 1$. This shows that each of these links are neither isotopic to the unlink, nor to each other. 

Further, we can generalize this for $m>3$ to show that 
\begin{center}
deg$_L(V(W_{n}(B_m))) = \displaystyle \frac{5-5m}{2} -(n-1)$ for $n \ge 1$ and $m >3$.
\end{center}

This shows that the links $W_n(B_m)$ are all non-isotopic links. 

 \vspace{3mm}

\section{General links}
\label{general}

In this section we will prove Theorem \ref{mainthm} for a general link. We first define what we mean by \textit{connect sum of links}.

\textbf{Definition 5.1} (Connect sum of links) 
Let $L = (l_1, \ldots , l_n)$ and $L' = (l'_1, \ldots , l'_n)$ be two $n$-component links in $S^3$. The links can be considered to be in $\R^3$.

\begin{wrapfigure}{r}{0.58\textwidth}
\centering
\includegraphics[width=0.85\linewidth]{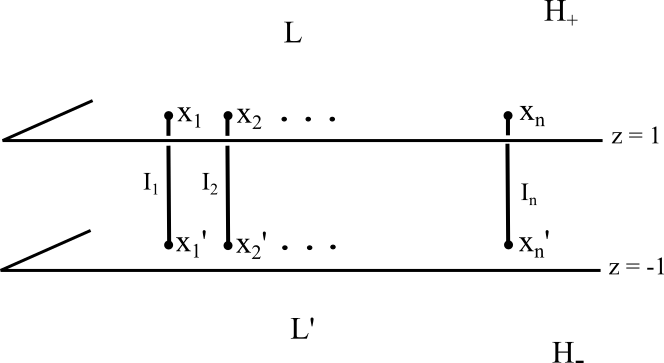}
\caption{}
\label{Connect1}
\end{wrapfigure}

 Let $H_+ = \{ (x, y,z) \in \R^3 : z > 1 \}$ and $H_- = \{ (x, y,z) \in \R^3 : z < -1 \}$ be two half-spaces in $\R^3$. Consider the set of $n$ points $\{ x_i \in \R^3 : x_i = (i, 0 , 1) , i = 1, 2, \ldots, n \}$ in $\partial H_+$ as shown in Figure \ref{Connect1}, and similarly another set $\{ x_i' \in \R^3 : x_i' = (i, 0 , -1) , i = 1, 2, \ldots, n \}$ in $\partial H_{-}$. Let $I_i$ be line segments in $\R^3$ joining $x_i$ and $x_i'$ for $i = 1, 2, \ldots , n$ (Figure \ref{Connect1}). 

\begin{wrapfigure}{r}{0.58\textwidth}
\centering
\includegraphics[scale=0.45]{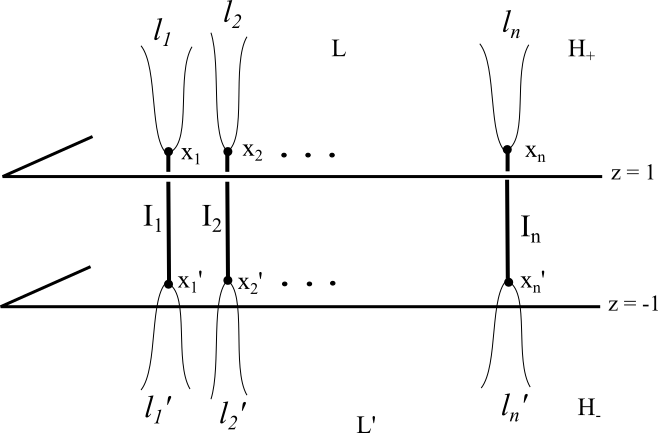}
\caption{}
\label{Connect2}
\end{wrapfigure}
Now we can assume that $L$ lies in $H_+$ and $L'$ in $H_-$. Let $\phi_+$ be an isotopy of $\R^3$ supported in the upper half space $\{ (x, y, z ) : z >0 \}$ such that $\phi_+ (L)$ lies in $\overline{H}_+$ and that each component $\phi (l_i)$ of $\phi (L)$ intersects $\partial H_+$ in exactly one point, $x_i$. Similarly, let $\phi_-$ be an isotopy of $\R^3$ supported in the lower half space such that $\phi_+ (L')$ lies in $\overline{ H}_-$ and that each component $\phi (l_i')$ intersects $\partial H_-$ in the point $x_i'$. See Figure \ref{Connect2}.

Since $\phi_+$ and $\phi_-$ are isotopies, we can call the new links $L$ and $L'$. A connect sum of each pair of components $l_i$ and $l_i'$ along the segment $I_i$ is formed in the following way: small unknotted arcs of $l_i$ and $l'_i$ are removed in small neighbourhoods of $x_i$ and $x_i'$. The endpoints of these arcs are connected by strands parallel to $I_i$ as shown. 

\begin{wrapfigure}{r}{0.58\textwidth}
\vspace{-8mm}
\centering
\includegraphics[scale=0.45]{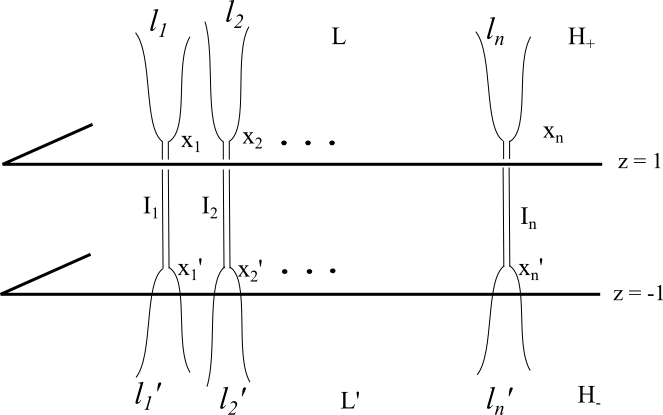}
\caption{}
\label{Connect3}
\vspace{-10mm}
\end{wrapfigure}

The resulting $n$-component link is called a connect sum of $L$ and $L'$ and denoted by $L \# L'$, as shown in Figure \ref{Connect3}.

If the links $L$ and $L'$ are oriented links, we can define $L\# L'$ to be an oriented link by having $h$ preserve the orientation of each component. 

\vspace{5mm}
Note that the definition of connect sum of links is not well-defined in general. It depends on  
\begin{itemize}
\item the ordering of the components,
\item the isotopies $\phi_+$ and $\phi_-$, 
\end{itemize}

\begin{figure}[h!]
\begin{subfigure}{0.45\textwidth}
\centering
\includegraphics[scale=0.7]{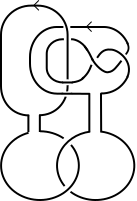}
\caption{}
\label{HCWa}
\end{subfigure}
\begin{subfigure}{0.5\textwidth}
\centering
\includegraphics[scale=0.55]{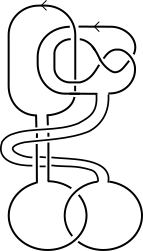}
\caption{}
\end{subfigure}
\caption{Connect sums of the Whitehead link and Hopf link}
\label{HCW}
\end{figure}

For example, Figure \ref{HCW} shows two different connect sums of the Whitehead link and the Hopf link. 

We first consider the case of 2-component links. Let $L$ be any given 2-component link in $S^3$. Let $H$ be the Hopf link. Then let $L \# H$ be a connect sum of $L$ and $H$. Recall that, by definition of connect sums, this involves making choices for the isotopies of half spaces contaning $L$ and $H$. We will call this choice the $type$ of the connect sum. 
Now replace one of the components of $H$ by its Whitehead $n$-double as defined in Def \ref{Whiteheaddouble} and call the resulting link $W_n(H)$. Construct a connect sum of $L$ and $W_n(H)$ using the same choice of isotopies, and hence of the same type. We will denote this particular choice of the connect sum by $L_n(H)$. 


\begin{wrapfigure}{r}{0.5\textwidth}
\centering
\includegraphics[scale=0.9]{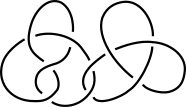}
\caption{2-component link $L$}
\label{L2}
\end{wrapfigure}

\begin{figure}[h!]
\begin{subfigure}{0.45\textwidth}
\centering
\includegraphics[scale=0.6]{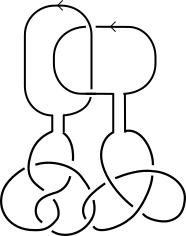}
\vspace{2mm}
\caption{$L \# H$}
\label{L(H)}
\end{subfigure}
\begin{subfigure}{0.5\textwidth}
\centering
\includegraphics[scale=0.6]{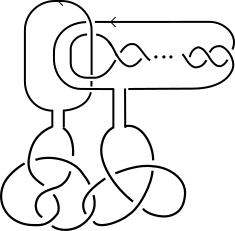}
\vspace{2mm}
\caption{$L_n(H)$}
\label{L_n(H)}
\end{subfigure}
\caption{}
\end{figure}

For example, let $L$ be the link with the Figure-8 knot and trefoil as its components as shown in Figure \ref{L2}. Let $L \# H$ be its connect sum with the Hopf link as showin in Figure \ref{L(H)}. Then the connect sum of a Whitehead $n$-double of $H$ with $L$ corresponding to the fixed isotopies will be as shown in Figure \ref{L_n(H)}. This gives us a family of connect sums of the same $type$. 

Let $L_n(H)$ be the connect sum shown in Figure \ref{L_n(H)}. 

Notice that when $n$ is even a single (self) crossing change in the Whitehead twists will reduce the link $L_n(H)$ to the original link $L$, showing that it is link homotopic to $L$. Also, if one component from $W_n(H) $  is removed then we get an unknot, and hence removing one component from the connect sum will give us the corresponding component of the link $L$. 
We will show that $L_n(H)$ is not isotopic to $L$ for large enough even integers $n$. Before we do that let us consider the case of links with more than two components. 


Let $L$ be any given $m$-component link in $S^3$ for $m \ge 3$. Let $B_m$ be the Brunnian link with $m$ components, and $L \# B_m$ be a connect sum. Just like in the case of 2-component links, we fix the isotopies for this connect sum. Consider the Whitehead $n$-double of one of the components of $B_m$ and call the resulting link $W_n(B_m)$. Now consider the connect sum $L \# W_n(B_m)$ with the same choices of isotopies and call it $L_n(B_m)$.


\begin{figure}[h!]
\begin{subfigure}{0.45\textwidth}
\centering
\vspace{11mm}
\includegraphics[scale=0.75]{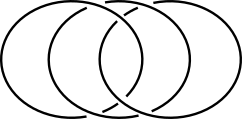}
\vspace{6mm}
\caption{3-component great circle link $L$}
\label{L3}
\end{subfigure}
\begin{subfigure}{0.5\textwidth}
\centering
\includegraphics[scale=0.5]{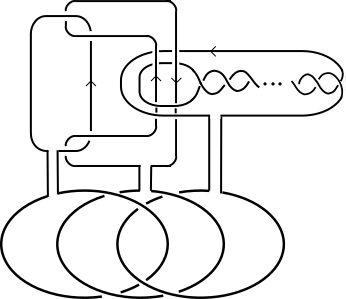}
\caption{Connect Sum $L_n = L \# W_n(B_3)$}
\label{L_n(B_3)}
\end{subfigure}
\caption{}
\end{figure}

For example let 
$L$ be the 3-component link with pairwise linking numbers $=1$  as shown in Figure \ref{L3} and let $L_n(B_3)$ be a family of connect sums of the same type as shown in Figure \ref{L_n(B_3)}. Just as in the case of 2-component links we can easily see that for $n$ even, $L_n(B_3)$ is link homotopic to $L$ and any proper sublink will be isotopic to the corresponding sublink of $L$.

 Now it remains to show that these links are indeed isotopically different from $L$. To do this we consider the case when the number of components $m = 3$.

\begin{figure}[h!]
\includegraphics[scale=0.5]{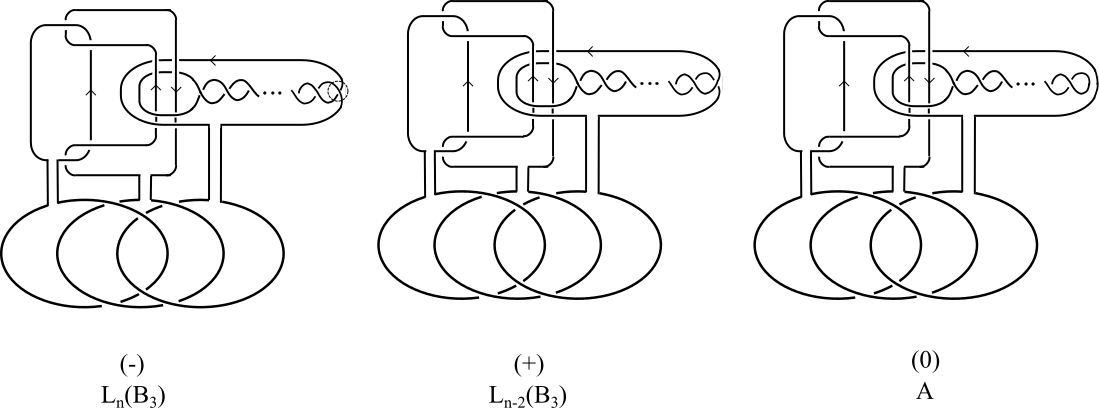}
\caption{}
\label{SkeinLnB3}
\end{figure}

We use Skein relations on the indicated crossing in Figure \ref{SkeinLnB3} to get 
\begin{equation}
t^{-1}V(L_n(B_3)) - tV(L_{n-2}(B_3)) + (t^{-1/2} - t^{1/2}) V(A) = 0,\label{A}\end{equation}
where $A$ is the link as shown in Figure \ref{SkeinLnB3}. We can compute the Jones polynomial of $A$ using Skein relations multiple times to get\\
$V(A) =  t^{23/2} - 3t^{21/2} + 5t^{19/2} - 2t^{17/2} - 11t^{15/2} + 15t^{13/2} - t^{6} - 17t^{11/2} + 8t^{5} + 20t^{9/2} - 15t^{4} - 28t^{7/2} + 12t^{3} + 21t^{5/2} - 7t^2 - 8t^{3/2} + 8t - 4t^{1/2} - 7 + 4t^{-1/2} + 4t^{-1} - 3t^{-3/2} - t^{-2} + t^{-5/2}.$

As before, consider the function 
$$\Phi: \Z [t^{-1/2}, t^{1/2}] \rightarrow \Z [t^{-1/2}, t^{1/2}] $$ 
$$\Phi(p) = t^{-2} p + t^{-1} (t^{-1/2} - t^{1/2}) V(A),$$ where $V(A)$ is the Jones polynomial of the link $A$ as given in equation (\ref{A}). 
Then $V(L_n(B_3)) = \Phi(V(L_{n-2}(B_3)))$, for every $n \ge 2$. Hence we only need to analyse the function $\Phi$ in order to understand the Jones polynomials of the links $L_n$.
Notice that the degree of the lowest degree term of $\Phi(p)$ can be given by
\begin{center} deg$_L(\Phi(p)) = $ min \{deg$_L (p) - 2, $ deg$_L (V(A)) - 3/2 $\}. \end{center}

We know from our computations that deg$_L (V(A)) = -5/2$ and hence, \begin{center} deg$_L(\Phi(p)) = $ min \{deg$_L (p) - 2,\ -4 $\}.\end{center}

For $n=2$, we have that $V(L_2(B_3)) = \Phi(V(L_0(B_3)))$. But $L_0(B_3)$ is nothing but the link $L$. We can compute the Jones polynomial of $L$, which comes out to be $2t^6 +t^4 - t^2$.

Hence deg$_L(V(L_0(B_3))) = 2$, and deg$_L(V(L_2(B_3))) =$ deg$_L(\Phi(V(L))) = -4$.
Continuing in a similar fashion, we get that deg$_L(V(L_4(B_3))) = -6,$ deg$_L(V(L_6(B_3))) = -8$ and so on. 

We can generalize this to obtain that deg$_L(V(L_{n}(B_3))) = -n-2$ for $n\ge 2$, $n$ even. This shows that each of the links $L_n$ are neither isotopic to $L$, nor to each other. 

Notice that it is crucial that we chose the connect sums $L_n(B_m)$ to be of the same type, otherwise the link $A$ might change at every stage of the iteration giving us different Jones polynomials $V(A)$. 

In general, the link $A$ depends on the link $L$ and the type of connect sum, and may have a different Jones polynomial. However, we still have that 
deg$_L(V(L_n(B_3))) = $ min \{deg$_L (V(L_{n-2}(B_3))) - 2, $ deg$_L (V(A)) - 3/2 $\} for $n \ge 2$. Hence, we will get a positive integer $N$ such that for every even $n > N$ 
deg$_L(V(L_n(B_3))) = $ deg$_L (V(L_{n-2}(B_3))) - 2$, showing that eventually all the Jones polynomials $V(L_n(B_3))$ will be distinct and hence the links will be non-isotopic.

\section{Milnor invariants}
\label{milnor}

In this section we give a proof of Theorem \ref{milnorinv}. Let $L$ be an $m$-component link in $S^3$. Let $F(L)$ be the fundamental group, $\pi_1(S^3  - L)$, and $F_q(L)$ be its $q$-th lower central subgroup. We now define Milnor's $\mu$-invariants \cite{Mil2}. 


We have a presentation of the group $F(L)/F_q(L)$ with $m$ generators, $\alpha_1,\ldots, \alpha_m$, where $\alpha_i$ represents the $i$-th meridian of the link $L$. For each $i = 1, \ldots, m$, the longitude $l_i$ of the $i$-th component of $L$ is represented by a word in $\alpha_1, \ldots, \alpha_m$. 

On substituting $\alpha_j = 1 + \kappa_j$, and $\alpha_j^{-1} = 1 - \kappa_j + \kappa_j^2 -\kappa_j^3 +\ldots $ in the expression for $l_i$, we can represent $l_i$ by a formal power series in $\kappa_1, \ldots, \kappa_m$. Let $\mu (j_1 \ldots j_s i)$ denote the coefficient of $\kappa_{j_1} \ldots \kappa_{j_s} $ in this series, for every $s\ge 1$. 

Let $\Delta (i_1 \ldots i_r)$ denote the greatest common divisor of $\mu (j_1 \ldots j_s )$, where $j_1 \ldots j_s (2 \le s < r)$ ranges over all sequences obtained by cancelling at least one of the indices $ i_1 \ldots i_r$, and permuting the remaining indices cyclically. 
And now $\overline{\mu}(i_1 \ldots i_r)$ is defined as the residue class of $\mu(i_1 \ldots i_r)$ modulo $\Delta (i_1 \ldots i_r)$. 
These are Milnor's isotopy invariants of links. 


Consider the link obtained from the Brunnian link $B_3$ by taking the Whitehead double twice, on the first and last component (Figure \ref{W_{n,m}^{1,3}(B_3)}). We denote this link by $W_{n,m}^{1,3}(B_3)$, where $k$ and $n$ are the number of twists on the first and third components respectively. 

\begin{figure}[h]
\includegraphics[scale=0.7]{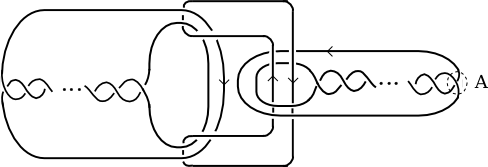}
\caption{$W_{k,n}^{1,3}(B_3)$}
\label{W_{n,m}^{1,3}(B_3)}
\end{figure}
By Corollary 4.2 in \cite{MY}, we know that all the Milnor invariants, $\overline \mu$, of $W_{k,n}^{1,3}(B_3)$ are 0. Hence these links satisfy property (4) of Theorem \ref{milnorinv}. However, we can no longer use Milnor invariants to show that the links are not isotopic to unlink. But we can use the Jones polynomials to show existence of links that satisfy all four properties.

It is clear that the links will satisfy properties (1) and (3). We show as before that the Jones polynomial will come out to be non-trivial showing that the links satisfy (2). 

On considering the Skein relation at the crossing (A) (Figure \ref{W_{n,m}^{1,3}(B_3)}) we get the Skein relation
$$t^{-1}V(W_{k,n-2}(B_3)) - tV(W_{k,n}(B_3)) + (t^{-1/2} - t^{1/2}) V(L_{k,3}) = 0,$$ 
where $L_{k,3}$ is the link in Figure \ref{L_{m,3}}.

\begin{figure}[h]
\includegraphics[scale=0.7]{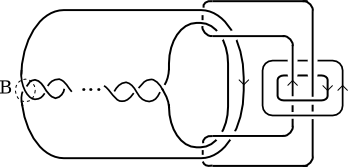}
\caption{$L_{k,3}$}
\label{L_{m,3}}
\end{figure}

If both $W_{k,n-2}(B_3)$ and $W_{k,n}(B_3)$ are unlinks then we will have $V(L_{k,3}) = (-t^{-1/2} - t^{1/2})^3.$

Now consider the Skein relation at crossing (B) in Figure \ref{L_{m,3}}. We get the following relation.
$$t^{-1}V(L_{k-2,3}) - tV(L_{k,3}) + (t^{-1/2} - t^{1/2}) V(L'_{3}) = 0,$$ 
where $L'_3$ is the following link.

\begin{figure}[h]
\includegraphics[scale=0.7]{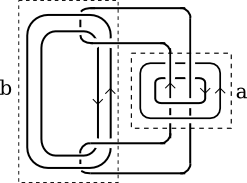}
\caption{$L'_3$}
\label{L'_3}
\end{figure}

We compute $V(L'_3)$ in similar fashion as we computed $V(L_3)$ in Section \ref{m=3}. We first apply Skein relations at the crossings shown in the box (a) in Figure \ref{L'_3}, and then on the crossings in box (b). 

This gives us the Jones polynomial $$V(L'_3) = - t^6 + t^5 + t^3 + 2t^2 + 2t + 6 + 2t^{-1} + 2t^{-2} + t^{-3} + t^{-5} - t^{-6}.$$

This shows that either $V(L_{k-2,3})$ or $V(L_{k,3})$ is different from $(-t^{-1/2} - t^{1/2})^3$. This in turn shows that at least one of $W_{k,n}(B_3)$, $W_{k,n-2}(B_3)$, $W_{k-2,n-2}(B_3)$ and $W_{k-2,n-2}(B_3)$ is isotopically non-trivial, showing the existence of links satisfying $(1) - (4)$ in Theorem \ref{milnorinv}. 

We can extend these computations to links $L'_m$ for $m>3$ to get the recurrence relation 
$$V(L'_m) = V(L'_{m-1}) (t^{5/2} -t^{3/2} -t^{-3/2} +t^{-5/2}) - (-t^{1/2}-t^{-1/2})^m (t^{5/2} -t^{3/2} -  t^{1/2} + t^{-1/2} -t^{-3/2} +t^{-5/2}).$$

Analysis similar to that in the previous section shows that $V(L'_m)$ is not equal to $(-t^{-1/2} - t^{1/2})^m$ (the degree of $V(L_m')$ being strictly larger than $m$), which in turn proves Theorem \ref{milnorinv} for $m$-component links with $m>3$.


\section{Concluding Remarks}
\subsection{Other invariants} 
Consider the Seifert matrix of the link $W_n^-(B_3)$, the Whitehead double of $B_3$ with $n$ negative twists, where $n\ge 4$ is even. 
\[
S_n = \left(
\begin{array}{cccc|c}
   & \hspace{5mm} \raisebox{-25pt}{{\huge\mbox{{$S$}}}} & & & 0 \\[-5ex]
  & & & & \vdots \\
  & & & & 0 \\ 
  & & & & 0 \\
  & & & & 1\\
  & & & & -1\\ \hline
  0 \cdots 0 & -1 & 0 & -1 & \frac{n-4}{2} 
\end{array}
\right)
\]
where $S$ is nothing but the Seifert matrix for the link $L_3$ from Figure \ref{Boxed L_3}.

On looking at the form of the matrix, one is tempted to compare the signatures of the two links. To do this, consider the symmetric bilinear form given by $A_n = S_n+S_n^T$.
The \textit{signature} of the link $W_n^-(B_3)$, $\sigma (W_n^-(B_3))$, is by definition, the signature of $A_n$. 
\[
A_n = \left(
\begin{array}{cccc|c}
 &  \hspace{5mm} \raisebox{-25pt}{{\huge\mbox{{$A$}}}} & & & 0 \\[-5ex]
  & & & & \vdots \\
  & & & & 0 \\ 
  & & & & -1 \\
  & & & & 1\\
  & &  & & -2\\ \hline
  0 \cdots 0 & -1 & 1 & -2 & n-4
\end{array}
\right),
\]
where $A$ gives the bilinear form corresponding to the link $L_3$. (Detailed discussion on Seifert form and signature of links can be found in \cite{Lic}, Chapters 6 and 8.)

Suppose $A$ defines a bilinear form on the vector space $V$, and $A_n$ on $V \oplus \R$. Let $w_n$ be the eigenvector for the largest eigenvalue, $\lambda_n$, of $A_n$, and $H_n$ be the subspace orthogonal to $w_n$. The form of the matrix $A_n$ allows us to conclude that $w_n \rightarrow (0, \ldots, 0, 1)$ and $\lambda_n \rightarrow \infty$, as $ n \rightarrow \infty$. This in turn gives us that subspace $H_n$ approaches $V$ as $n\rightarrow \infty$. Hence, the form $A_n|_{H_n}$ approaches $A_n|_{V}$. However, $A_n|_{V}$ is same as the form $A$.

This shows that the eigenvalues of $A_n|_{H_n}$ must also approach those of $A$. Thus, \textit{as long as none of the eigenvalues of $A$ are zero}, for a large enough value of $n$, we would obtain $\sigma(A_n|_{H_n}) = \sigma(A)$, which would force $\sigma(A_n) = \sigma(A) + 1$. Alternatively, if $W_n^+(B_3)$ is a link having $n$ positive twists with $n \ge 4 $ being even, then a similar analysis shows that for large $n$, $\sigma (A_n ) = \sigma(A) - 1.$

Thus, for large enough $n$ even, $\sigma(W_n^+(B_3)) = \sigma(L_3) - 1$ 
and $\sigma(W_n^-(B_3)) = \sigma(L_3) + 1$. Hence at least one of the links, $W_n^+(B_3)$ or $W_n^-(B_3)$, must have signature different from that of the unlink, proving the existence of a link as in Theorem \ref{mainthm} for $m=3$. The same method would work for links $W_n(B_m)$ for $m >3$ as well. \\

Unfortunately, the following observation shows that det$(A)=0$, so $A$ will have zero as an eigenvalue, showing that the above method will not work for these examples of links.

A \textit{coloring} (modulo $n$) of a link $L$ is a function $c$ from the set of arcs in its link diagram to $\Z /n\Z$ if at every crossing we have $2c(x_i) \equiv c(x_j) + c(x_k)$ mod $n$, where $x_i, x_j, x_k$ are arcs at a crossing as shown.

\begin{figure}[h]
\includegraphics[scale=0.7]{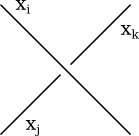}
\caption{}
\end{figure}

A system of these equations gives us the coloring matrix. Deleting any one column and one row yields a new matrix. We know that det$(A)$ is equal to the determinant of this new matrix for link $L_m$. And this determinant is $0$ if and only if the link has a coloring mod $n$ for every $n \in \Z$, and in particular, if it has a coloring over $\Z$ (\cite{L}, Chapter 3).

For every $m \ge 3$, we have the following coloring of $L_m$ over $\Z$, which forces det$(A) = 0$. 

\begin{figure}[h]
\includegraphics[scale=0.7]{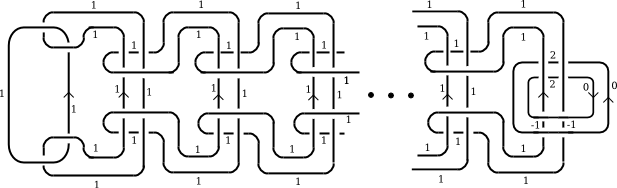}
\caption{}
\end{figure}

\subsection{3-manifold theory} 
 As per the comments by an anonymous referee, one can use 3-manifold theory to show that the Whitehead double of (any non-split) link is isotopically non-trivial, and in fact, non-split. The complement of a Brunnian link in $S^3$ gives an irreducible 3-manifold with incompressible boundary, and so does the complement of the Whitehead $n$-twists in the solid torus. 

Let $\mathcal W_n$ be the Whitehead $n$ twists in the solid torus $D^2 \times S^1$ (Figure \ref{Whitehead double}). Let $M_1 \defeq S^3 \setminus N (B_m)$ and $M_2 \defeq (D^2 \times S^1) \setminus N (\mathcal W_n)$, where $N(B_m)$ and $N(\mathcal W_n)$ denote the tubular neighbourhoods. Let $T^2$ be the torus boundary of the tubular neighborhood of the $m$-th component of $B_m$, $N (B_m^m)$. Then we have $S^3 \setminus N(W_n(B_m)) =  M_1 \cup_{T^2 } M_2$. Since $M_1$ and $M_2$ are irreducible manifolds, gluing them along their incompressible boundary gives us another irreducible manifold. This shows that $W_n(B_m)$ cannot be a trivial link, since otherwise we get a 2-sphere that separates $S^3 \setminus N(W_n(B_m))$. 

To see that gluing $M_1$ and $M_2$ along their incompressible boundary gives us another irreducible manifold, say $M$, let us first assume that $M$ is not irreducible. Then there must be a copy of $\alpha : S^2 \hookrightarrow M
$ which separates $M$. Also, $\alpha(S^2)$ must intersect $T^2$ , since otherwise it would lie completely in either $M_1$ or $M_2$, which are irreducible. We may assume that 
$\alpha$ is transverse to $T^2$. Then $\alpha (S^2) \cap T^2$ is a codimension 1 submanifold in $T^2$. Then $\alpha ^{-1} (T^2) \hookrightarrow S^2$ also has codimension 1, and hence must be a union of copies of $S^1$ in $S^2$. 

Consider a copy of $S^1$, say $S^1_\beta$, in $\alpha(S^2) \cap T^2$ that bounds a disk in $M_1$. Since $T^2$ is incompressible in $M$, it must bound a disk in $T^2$. This gives a copy of $S^2 $ in $M_1$. Being the complement of the link $B_m$, $M_1$ is aspherical. Hence the 2-sphere must bound a 2-ball in $M_1$. Then the map $\alpha$ can be perturbed so that $\alpha(S^2) \cap T^2$ no longer contains $S^1_\beta$. Doing this for all the copies of $S^1$ in $S^2 \cap T^2$ that bound a disk in $M_1$ will eventually show that there is a copy of $S^2$ in $M_2$ that separates $M_2$ which contradicts the irreducibility of $M_2$. 
 
We note here that a similar 3-manifold theory technique does not help in proving the result for a general link $L$. 
In the case of a general link, we need to show that the connect sum is not isotopic to the link $L$, which is different from showing that there is a separating 2-sphere in $S^3$. Further, the decomposition of $S^3 \setminus N(L\# W_n(B_m))$ is a bit more complicated than the unlink case and similar techniques don't work anymore.

We further note that the Jones polynomial methods help us to show that all the links $L \# W_n(B_m)$ are isotopicallly distinct giving us an infinite family of such links.

\end{document}